\documentclass[12pt]{article}

\usepackage{amscd,amsmath, amssymb, fancyhdr, mathbbol, dutchcal, url}
\usepackage[backref=page]{hyperref}
\renewcommand*{\backrefalt}[4]{%
	\ifcase #1 (Not cited.)%
	\or        (Cited on page~#2.)%
	\else      (Cited on pages~#2.)%
	\fi}

\hypersetup{
	colorlinks   = true,
	citecolor    = magenta,
	linkcolor    =blue,
	urlcolor     =magenta	
}

\numberwithin{equation}{section}

\newcommand{\version}{version 4,\ \ June 21, 2021}

\def\eqref#1{(\ref{#1})}

\newcommand{\arrow}{{\:\longrightarrow\:}}

\def\C{{\Bbb C}}
\newcommand{\R}{{\Bbb R}}

\def\1{\sqrt{-1}\:}
\newcommand{\restrict}[1]{{\left|_{{\phantom{|}\!\!}_{#1}}\right.}}
\newcommand{\cntrct}                
{\hspace{2pt}\raisebox{1pt}{\text{$\lrcorner$}}\hspace{2pt}}

\renewcommand{\phi}{\varphi}
\renewcommand{\epsilon}{\varepsilon}

\newcommand{\Vol}{\operatorname{Vol}}

\newcommand{\Sym}{\operatorname{Sym}}

\newcommand{\const}{\operatorname{\sf const}}

\newcommand{\Disc}{\operatorname{Disc}}

\newcounter{Mycounter}[section]
\newcounter{lemma}[section]
\renewcommand{\thelemma}{{Lemma \thesection.\arabic{lemma}}}
\newcommand{\lemma}{%
    \setcounter{lemma}{\value{Mycounter}}
    \refstepcounter{lemma}
    \stepcounter{Mycounter}
    {\noindent \bf \thelemma:\ }}

\newcounter{claim}[section]
\renewcommand{\theclaim}{{Claim \thesection.\arabic{claim}}}
\newcommand{\claim}{%
    \setcounter{claim}{\value{Mycounter}}
    \refstepcounter{claim}
    \stepcounter{Mycounter}
    {\noindent \bf \theclaim:\ }}

\newcounter{sublemma}[section]
\newcounter{corollary}[section]
\newcounter{theorem}[section]
\renewcommand{\thetheorem}{{Theorem \thesection.\arabic{theorem}}}
\newcommand{\theorem}{%
    \setcounter{theorem}{\value{Mycounter}}
    \refstepcounter{theorem}
    \stepcounter{Mycounter}
    {\noindent \bf \thetheorem:\ }}

\newcounter{conjecture}[section]
\newcounter{proposition}[section]
\newcounter{definition}[section]
\renewcommand{\thedefinition}
      {{Definition~\thesection.\arabic{definition}}}
\newcommand{\definition}{%
    \setcounter{definition}{\value{Mycounter}}
    \refstepcounter{definition}
    \stepcounter{Mycounter}
    {\noindent \bf \thedefinition:\ }}

\newcounter{example}[section]
\newcounter{remark}[section]
\renewcommand{\theremark}{{Remark \thesection.\arabic{remark}}}
\newcommand{\remark}{%
    \setcounter{remark}{\value{Mycounter}}
    \refstepcounter{remark}
    \stepcounter{Mycounter}
    {\noindent \bf \theremark:\ }}

\newcounter{problem}[section]
\newcounter{question}[section]
\newcommand{\proof}{\noindent{\bf Proof:\ }}

\makeatletter

\@addtoreset{equation}{section} \@addtoreset{footnote}{section}
\makeatother

\def\blacksquare{\hbox{\vrule width 5pt height 5pt depth 0pt}}
\def\endproof{\blacksquare}

\begin{document}
\begin{center}
{\LARGE\bf
Holomorphic Lagrangian subvarieties in holomorphic symplectic manifolds
with Lagrangian fibrations and special K\"ahler geometry\\[4mm]
}

Ljudmila Kamenova\footnote{Partially supported 
by a grant from the Simons Foundation/SFARI (522730, LK)}, 
Misha Verbitsky\footnote{Partially supported 
by  the  Russian Academic Excellence Project '5-100',
FAPERJ E-26/202.912/2018 and CNPq - Process 313608/2017-2.}

\end{center}

{\small \hspace{0.10\linewidth}
\begin{minipage}[t]{0.85\linewidth}
{\bf Abstract.} 
Let $M$ be a holomorphic symplectic K\"ahler manifold 
equipped with a Lagrangian fibration $\pi$ with compact fibers. 
The base of this manifold is equipped with a {\bf special K\"ahler
structure}, that is, a K\"ahler structure $(I, g, \omega)$ and a symplectic
flat connection $\nabla$ such that the metric $g$ is locally the Hessian of
a function. 
We prove that any Lagrangian subvariety $Z\subset M$
which intersects smooth fibers of $\pi$ and smoothly projects to $\pi(Z)$ 
is a torus fibration over its image $\pi(Z)$ in $B$, and this image is
also special K\"ahler. This answers a question of N. Hitchin related
to Kapustin-Witten BBB/BAA duality.
\end{minipage}
}


\section{Introduction} 


The present paper is motivated by the observations made by N. Hitchin
\cite{_Hitchin:talk_} who worked on the Kapustin-Witten
version of the geometric Langlands correspondence,
interpreted as Montonen-Olive generalization of 
electric-magnetic duality. This theory originates in
1977, when P. Goddard, J. Nuyts and D. Olive
discovered that magnetic sources in gauge theory with gauge group $G$
are classified by irreducible representations of the Langlands
dual group ${}^LG$ (\cite{_Goddard_Nuyts_Olive_}). 
Then C. Montonen and D. Olive conjectured that the Yang-Mills theories
with the gauge groups $G$ and ${}^LG$ are isomorphic on
the quantum level. The Montonen-Olive duality can be
regarded as a quantum field generalization of the usual
electric-magnetic duality.

M. Atiyah had suggested that the Montonen-Olive
conjecture (\cite{_Montonen_Olive:monopoles_})
might be related to the Langlands duality, but 
it took many years until 2006, when A. Kapustin and E. Witten 
explained this conjectural relation.

In their celebrated paper \cite{Kapustin_Witten}, Kapustin and Witten produced
a rich dictionary of the correspondence between the geometric Langlands program 
and S-duality in the four-dimensional N=4 gauge theory.
This approach is based on the comparison between two
Hitchin systems (the spaces of Higgs bundles on a curve)
with values in Langlands dual groups. Both of these
Hitchin systems are equipped with a Lagrangian fibration.
Reminiscent of the Strominger-Yau-Zaslow interpretation
of the Mirror Symmetry, the Langlands duality
is interpreted as a correspondence between certain
categories on these two spaces,
associated with the duality of their fibers.
For a less technical survey of the Kapustin-Witten
program, see \cite{_Kapustin:survey_}.

The Kapustin-Witten interpretation
of Montonen-Olive/geometric Langlands duality can be understood
as SYZ Mirror symmetry on the Hitchin space,
but it is firmly based on the hyperk\"ahler geometry of the Hitchin space.
In place of the Fukaya category on the symplectic
side of Mirror Symmetry, one has
a category associated with the holomorphic
Lagrangian subvarieties (BAA, ABA and AAB branes).
In place of the derived category of coherent sheaves on the complex side
of Mirror Symmetry one has a category which has pairs (trianalytic subvariety,
hyperholomorphic bundle on it) as objects; these are called BBB branes. 
Since the fiberwise duality should somehow exchange 
these two categories, Hitchin argued, the fibers 
of the BAA brane under the Hitchin fibration map
should be tori, and its image 
should retain the special K\"ahler structure which
exists on the base of the Hitchin fibration. 
We define all these notions and state this result 
rigorously in Section \ref{_Special_K_Section_}.

Hitchin stated his theorem in bigger generality
than required by the Kapustin-Witten theory:
he expected it to be true for any hyperk\"ahler
manifold equipped (such as the Hitchin system) with a $\C^*$-action
rotating the complex structures within the twistor family.
We prove the same result without a $\C^*$-action. Our main theorem is the
following. 

\hfill

\theorem (\ref{_subva_main_Theorem_})
Let $(M, \Omega)$ be a holomorphic symplectic K\"ahler manifold, 
and let $\pi:\; M \arrow B$ be a proper Lagrangian fibration. Consider
an irreducible Lagrangian subvariety $Z\subset M$ such that $\pi(Z)$ does not
lie in the discriminant locus $D$ of $\pi$.  Then for any smooth point
$x\in \pi(Z)\backslash D$ which is a regular value of $\pi:\; Z\arrow \pi(Z)$, 
the fiber $\pi^{-1}(x)\cap Z$ is a union of translation 
equivalent subtori in the complex torus $\pi^{-1}(x)$, 
and the regular part of $\pi(Z)$ is a special K\"ahler submanifold in 
$B\backslash D$.


\section{Special K\"ahler manifolds} 
\label{_Special_K_Section_}


\subsection{Special K\"ahler manifolds and Hessian manifolds}

Special K\"ahler manifolds first appeared in physics,
\cite{_deWit_vanProyen_}, \cite{_deWit_Lauwers_vanProyen_},
as allowed targets for the scalars of the vector multiplets
of field theories with $N=2$ supersymmetry
on a 4-dimensional Minkowski space-time. Originally
they came in two flavours, the affine special K\"ahler
manifolds associated with rigid supersymmetry, and
projective special K\"ahler manifolds associated with
the local supersymmetry. In the present paper we are
interested only in the affine version.

The first comprehensive mathematical exposition of this theory
is due to Dan Freed, \cite{_Freed:SK_}.
After this geometric structure was presented to the general
mathematical readership, special K\"ahler manifolds became
prominent in differential geometry. In \cite{_Baues:Cortes_}, 
Baues and Cort\'es have
shown that special K\"ahler manifolds can be
interpreted as ``affine hyperspheres''.
This classical concept, going back to the work of Blaschke
in affine geometry, is described by solutions of real
Monge-Amp\`ere equation. This interpretation leads
to a classification of special K\"ahler manifolds.
For more details on the differential geometry
of special K\"ahler manifolds, the reader is
directed to the survey \cite{_Cortes:survey_}. 

\hfill

\definition
A {\bf special complex manifold}
is a complex manifold $(M,I)$ equipped with a 
flat, torsion-free connection $\nabla$ such that
the tensor $\nabla(I)\in \Lambda^1(M)\otimes\Lambda^1(M)\otimes TM$
is symmetric in the first two variables.
A special complex manifold is {\bf special K\"ahler}
if it is equipped with a K\"ahler form $\omega$
which satisfies $\nabla(\omega)=0$.

\hfill

Let $(M, I, \nabla, g, \omega)$ be a special K\"ahler manifold.
Since $\nabla(\omega)=0$, and $\nabla(I)$ is symmetric
in the first two variables, the tensor 
\begin{equation} \label{_deriva_g_Equation_}
\nabla(g)= \nabla(I \circ \omega)\in \Lambda^1(M)\otimes\Lambda^1(M)\otimes \Lambda^1(M)
\end{equation}
is symmetric in the first two variables. This tensor is symmetric
in the last two variables, because $g$ is symmetric. Therefore,  
$\nabla g$ is a symmetric 3-tensor.

\hfill

\definition
Let $(M, \nabla)$ be a manifold equipped with a flat torsion-free connection,
and $g$ a Riemannian metric. It is called {\bf Hessian} if $\nabla(g)$ is
symmetric in all 3 variables.

\hfill

\remark
It is not hard to see that the Riemannian metric $g$ on $(M,\nabla)$ 
is Hessian if and only if
$g$ is locally the Hessian of a function, which is called the {\bf potential} 
of the Hessian metric. A priori the potential exists only locally,
but when $M$ is simply connected, it can be defined globally on $M$.

\hfill

This construction is due to N. Hitchin, \cite{_Hitchin:SpLag_},
who exhibited the special K\"ahler
structure on the moduli space of holomorphic Lagrangian
subvarieties in a hyperk\"ahler manifold, and
exhibited many interesting differential-geometric
properties of special K\"ahler manifolds.

\hfill

\remark
Let $(M, I, \nabla, g, \omega)$ be a special K\"ahler manifold, $\Vol M$
the Riemannian volume form,
and $f$ the potential of its Hessian metric. Since $\Vol M=\omega^n$,
and $\nabla (\omega)=0$, the function $f$ is a solution of the real 
Monge-Amp\`ere equation $\det \frac{d^2f}{dx_idx_j}=\const$. In the paper
\cite{_Cheng-Yau:MA_}, Cheng and Yau studied Hessian manifolds 
with a prescribed Riemannian volume form, and proved an analogue 
of Calabi-Yau's theorem for such manifolds.

\hfill

\claim
Let $(M, I, \nabla, g, \omega)$ 
be a K\"ahler manifold equipped with a flat connection $\nabla$ which
satisfies $\nabla(\omega)=0$. Then  $(M, I, \nabla, g, \omega)$ 
is special K\"ahler if and only if the metric $g$ is Hessian.

\hfill

\proof Follows immediately from \eqref{_deriva_g_Equation_}. \endproof

\subsection{Special K\"ahler structure on the base of a complex Lagrangian fibration}

Special K\"ahler manifolds naturally occur in many situations associated with
the geometry of Calabi-Yau and hyperk\"ahler subvarieties. For the present 
paper, the following construction is most significant.

\hfill

\definition
Let $(M, \Omega)$ be a holomorphic symplectic manifold.
A {\bf Lagrangian subvariety} of $M$ is a subvariety such that its
smooth part is a Lagrangian submanifold in $M$. 
A (holomorphic) {\bf Lagrangian fibration} on $M$ is a proper holomorphic map
$\pi:\; M \arrow B$ with general fibers being Lagrangian submanifolds in $(M,\Omega)$.

\hfill

The following claim is well-known in classical mechanics.

\hfill

\claim
A smooth fiber of a holomorphic Lagrangian fibration is always a torus.

\hfill

\proof
For any fibration $\pi:\; M \arrow B$, any smooth fiber $F$ has
trivial normal bundle $NF$. However, $NF$ is dual to the tangent bundle
$TF$ whenever $\pi$ is a Lagrangian fibration. Therefore, 
the bundle $TF$ is also trivial. For any function on $B$,
its Hamiltonian gives a section of $TF$. Choose a collection
of holomorphic functions such that their Hamiltonians give
a basis in $TF$. Since these Hamiltonians commute, the corresponding
vector fields in $TF$ also commute. This gives a locally free action of an abelian
Lie group on $F$, and therefore $F$ is a quotient of an abelian
group by a lattice.
\endproof

\hfill

\definition
Let $\pi:\; M \arrow B$ be a proper fibration.
Consider the first derived direct image
$R^1\pi_*(\R_M)$, where $\R_M$ is the trivial sheaf on $M$.
This is a constructible sheaf; at 
any point $x\in B$, the fiber of 
$R^1\pi_*(\R_M)$ is equal to the first cohomology of the fiber $\pi^{-1}(x)$. 
Outside of singularities of $\pi$,
this sheaf is locally constant. The
flat connection on the corresponding vector bunlde is called 
the {\bf Gauss-Manin connection}. This connection is
defined in the complement to the set $\Disc(\pi)$ 
of all critical values of $\pi$; this set is
called the {\bf discriminant locus} of $\pi$. 

\hfill

\definition \label{_connection_base_lagr_Remark_}
Let $\pi:\; M \arrow B$ be a Lagrangian fibration,
and let $F$ be the fiber over $x\in B$. Then $\pi^* TB = NF= T^*F$. 
Identifying $H^0(NF)= T_x B$ with $H^0(T^*F)= H^1(F, \R)$, 
we obtain an identification of $TB$ and the bundle $R^1 \pi_* (\R_M)$ 
of the first cohomology constructed above. Therefore, 
$TB$ is equipped with a natural flat connection, 
also called the {\bf Gauss-Manin connection}.

\hfill

\remark\label{_2-form_base_lagr_Remark_}
 Let $\pi:\; M \arrow B$ be a holomorphic Lagrangian fibration. 
A K\"ahler form $\omega$ on $M$ restricted to 
a smooth fiber $F$ of $\pi$ defines a cohomology class
$[\omega]\in H^2(F)$. Since $F$ is a torus, we can
consider $[\omega]$ as a 2-form on $R^1\pi_*(\R_M)=TB$.
This form is clearly parallel under the Gauss-Manin connection.
Abusing the notation, we denote this 2-form by the same letter $\omega$.

\hfill

\theorem (\cite[Theorem 3.4]{_Freed:SK_}, \cite[Theorem 3]{_Hitchin:SpLag_}) 
Let $\pi:\; M \arrow B$ be a holomorphic Lagrangian fibration on
a K\"ahler holomorphic symplectic manifold
and let $B_0\subset B$ be the complement to the discriminant locus of $B$.
Consider the 2-form $\omega$ on $B$ constructed in 
\ref {_2-form_base_lagr_Remark_}, and the Gauss-Manin connection
$\nabla$  on $TB$ defined in \ref{_connection_base_lagr_Remark_}.
Then $(B, \nabla, \omega)$ is a special K\"ahler manifold.


\section{Special K\"ahler geometry and holomorphic Lagrangian subvarieties} 


\subsection{Holomorphic Lagrangian subvarieties: main theorem}

Recall that {\bf projective special K\"ahler manifold}
(\cite{_Mantegazza:spKahler_}) is a special K\"ahler manifold
$(M,g,I, \omega)$  equipped with a vector field $v$ acting on 
$(M,g)$ by homotheties which preserve the complex structure,
such that the vector field $I(v)$ acts by isometries.

\hfill

In his talk \cite{_Hitchin:talk_} at the SCGP in October 2018, 
Nigel Hitchin stated the following theorem.

\hfill

\theorem (Hitchin, \cite{_Hitchin:talk_}) 
Let $\pi\; : M \arrow B$ be an algebraically integrable system
with a $\C^*$-action defining a projective special K\"ahler structure. 
Then any $\C^*$-invariant holomorphic Lagrangian submanifold
has an open set with the structure of a fibration over
a projective special K\"ahler submanifold of $B$, and each fiber
is a disjoint union of translates of an abelian subvariety.

\hfill

Hitchin asked whether there is an analogue of his result in 
the affine special K\"aher setting. 
Here we prove it, and give examples of holomorphic Lagrangian submanifolds
projecting to special K\"ahler submanifolds.

\hfill

\theorem\label{_subva_main_Theorem_}
Let $(M, \Omega)$ be a holomorphic symplectic K\"ahler manifold, 
and let $\pi:\; M \arrow B$ be a proper Lagrangian fibration. Consider
an irreducible Lagrangian subvariety $Z\subset M$ such that $\pi(Z)$ does not
lie in the discriminant locus $D$ of $\pi$.  Then for any smooth point
$x\in \pi(Z)\backslash D$ which is a regular value of $\pi:\; Z \arrow \pi(Z)$,
the fiber $\pi^{-1}(x)\cap Z$ is a union of translation 
equivalent subtori in the complex torus $\pi^{-1}(x)$, 
and the smooth part of $\pi(Z)\backslash D$ is a special K\"ahler
submanifold in $B_0:=B\backslash D$.

\hfill

Before we prove \ref{_subva_main_Theorem_}, we state the 
following elementary linear-algebraic lemma.

\hfill

\lemma\label{_Lagra_linear_Lemma_}
Suppose $V\subset W\oplus W^*$ is a Lagrangian
vector subspace in $W\oplus W^*$ with standard
symplectic structure, and $\pi:\; W\oplus W^*\arrow W$
the projection. Then $\pi(V)^\bot= V\cap W^*$,
where $R^\bot\subset W^*$ denotes the annihilator of a subspace
$R\subset W$, and $W^*$ is considered as a subspace in $W\oplus W^*$.
\endproof

\hfill

{\bf Proof of \ref{_subva_main_Theorem_}. Step 1:}
Let $Z_x:= \pi^{-1}(x)\cap Z$, where $x\in \pi(Z)\backslash D$
is a regular value of $\pi:\; Z \arrow \pi(Z)$.
Denote by $T^\pi M$ the fiberwise tangent
bundle, and let $T^\pi_zM$ be its fiber over $z\in Z_x$.
The holomorphic symplectic form induces
non-degenerate pairing between $T^\pi_zM$ and $T_x B$.
\ref{_Lagra_linear_Lemma_}
applied to $V= T_z Z$ and $W\oplus W^*=T_z M$ implies 
$T_z Z_x=\pi(T_z Z)^\bot$.
Indeed, in this case $T_z Z_x=  V\cap W^*$ and
$\pi(T_z Z)=\pi(V)$. However,  $\dim \pi(Z)= \dim Z-\dim Z_x=\dim \pi(T_z Z)$, hence
in all points $x\in B\backslash D$ and all $z\in \pi^{-1}(x)$, one has
$T_x \pi(Z)= \pi(T_z Z)$. This gives
\begin{equation}\label{_orthogonal_complement_Lagra_Equation_} 
T_z Z_x= T_x(\pi(Z))^\bot.
\end{equation}

{\bf Step 2:}
From \eqref{_orthogonal_complement_Lagra_Equation_}, we obtain that
$\pi(T_z Z)^\bot$ is constant: for different $z, z'\in Z_x$, the spaces $T_z Z_x$ 
and $T_{z'} Z_x$ are obtained by a translation within the torus $\pi^{-1}(x)$.
In other words, the space $T_z Z_x$ is constant
in the standard coordinates on the torus, and
$Z_x\subset \pi^{-1}(x)$ is a union of
subtori which are translates of each other.

\hfill

{\bf Step 3:} Since $\pi(Z)\subset B$ is a complex subvariety,
in order to prove that it is special K\"ahler it suffices to show that
it is totally geodesic (that is, constant) with respect to the Gauss-Manin
connection $\nabla$ on $TB_0$.  However,
the connection $\nabla$ is identified with the
Gauss-Manin connection under the identification 
$TB_0= R^1\pi_*(\R_M)$, and it preserves any
sublattice in $T_xB_0=H^1(F, \R)$, where $F=\pi^{-1}(x)$.

Since $Z_x\subset \pi^{-1}(x)$ is a subtorus,
it corresponds to a sublattice $H_1(Z_x)\subset H_1(\pi^{-1}(x))$ 
in homology, and in a neighbourhood
 $U\ni x$, all fibers $Z_{x'}\subset \pi^{-1}(x')$
correspond to the same sublattice. Therefore,
its orthogonal complement  in $R^1\pi_*(\R_M)$ is constant.
However, by \eqref{_orthogonal_complement_Lagra_Equation_},
this orthogonal complement generates $T_x(\pi(Z))$.
This implies that $T \pi(Z)$ is constant with respect to
the Gauss-Manin connection $\nabla$ on $B_0$.
\endproof

\section{Examples} 

Many (or most) holomorphic Lagrangian tori in
hyperk\"ahler manifolds occur as fibers of Lagrangian fibrations.
Indeed, in \cite{_Hwang_Weiss_} it was shown that
any Lagrangian subtorus in a hyperk\"ahler manifold is a fiber
of a holomorphic Lagrangian fibration. However,
for any two distinct Lagrangian fibrations over
a maximal holonomy hyperk\"ahler manifold, the 
intersection index of their fibers is 
positive (\cite{_Kamenova_Lu_Verbitsky_}, second paragraph of the 
proof of Theorem 2.11). Therefore, any fiber of the first fibration
is projected to the base of the second one surjectively and 
finitely in the general point. In this case, 
\ref{_subva_main_Theorem_} is tautologically true, because
the fibers of $\pi\restrict Z$ are 0-dimensional, 
and its base coincides with $B$.

\hfill

It is much harder to find examples where the special
K\"ahler geometry of $\pi(Z)$ is non-trivial.
This is easy to explain. Indeed, $\pi(Z)$ gives
a flat submanifold in the special K\"ahler manifold $B_0=B\backslash D$.
Therefore, the tangent space $T_x \pi(Z)$ to any smooth point
is fixed by the monodromy of the Gauss-Manin connection on $B_0$.
However, the monodromy representation is quite often irreducible,
or has very few subrepresentaions. 

\hfill

The Hitchin system (moduli of Higgs bundles over a curve)
is equipped with a Lagrangian fibration (``Hitchin fibration''),
which has abelian varieties (Jacobians of the ``spectral curve'')
as its fibers. For some examples of the Hitchin
system, the corresponding monodromy representation 
was computed in \cite{_Baraglia_Schaposnik:monodro_}.
From this computations it follows that the monodromy representation
is reducible (\cite[Corollary 4.23]{_Baraglia_Schaposnik:monodro_}). 
This suggests that some interesting Lagrangian
subvarieties, not transversal to fibers of the Hitchin system,
might exist in this case. Two of the first papers to study the
monodromy for the Hitchin fibration are \cite{_Copeland_} and
\cite{_Schaposnik_}. 

\hfill

Holomorphic Lagrangian fibrations on a deformation
of the second Hilbert scheme of a K3 were studied by D. Markushevich and
by L. Kamenova (\cite{_Markushevich:gen2_}, \cite{Kamenova}). 
They split into two distinct cases. In the first case,
studied by D. Markushevich, the Abelian fibers
are Jacobians of smooth genus two curves. If the fibers of 
$\pi:\; M \arrow B$ have no elliptic curve then all 
Lagrangian subvarieties of $M$ would either project to $B$ 
surjectively or would lie in the fibers of $\pi$. 
In the second case, studied by L. Kamenova,
the fibers of $\pi$ are products of two elliptic curves, i.e., 
the fibers are Jacobians of singular genus two curves. 
This situation occurs, for example, when one takes the punctual 
Hilbert scheme of $2$ points on an elliptic K3 surface 
$S \arrow {\mathbb P}^1$. Then the Hilbert scheme $S^{[2]}$ is fibered over 
the base $({\mathbb P}^1)^{[2]} = {\mathbb P}^2$ with general fibers that are 
products of the fibers of the ellitic fibration $S \arrow {\mathbb P}^1$. 
As shown in \cite{Kamenova}, 
under some ``genericity'' hypotheses, all deformations
of the second Hilbert scheme of a K3, fibered with the fiber that is a 
product of two elliptic curves, are obtained in this way. 

\hfill

The main (and, so far, the only) non-trivial example of the 
geometric construction obtained in this paper is given  
by the Hilbert scheme of an elliptic K3 surface as follows. 
Let $\pi:\; M \arrow S=\C P^1$ be an elliptic fibration on a K3 surface.
Consider the corresponding fibration $\pi^{[n]}:\; M^{[n]}\arrow S^{[n]}=\C P^n$
on its Hilbert scheme. A {\bf multisection} of $\pi$ is
a curve which is transversal to the fibers of $\pi$; a multisection exists
if and only if $M$ is projective. Fix points $s_1, ..., s_k\in S$, and let
$C_{k+1}, ..., C_{n}\subset M$ be multisections.
Denote by $\hat L_k(s_1, ..., s_k, C_{k+1}, ..., C_n)\subset \Sym^n M$
the set of $n$-tuples of points $(e_1, ..., e_k, c_{k+1}, ..., c_n)\in \Sym^n M$,
such that $e_i\in \pi^{-1}(s_i)$ and $c_j \in C_j$. Since 
the holomorphic symplectic form on $\Sym^n M$ is locally a product
of the holomorphic symplectic form on $M$, and the curves
$\pi^{-1}(s_i)$ and $C_j\subset M$ are Lagrangian, the subvariety
\[ \hat L_k(s_1, ..., s_k, C_{k+1}, ..., C_n)\subset\Sym^n M\]is Lagrangian. 
Then its proper preimage
$L_k(s_1, ..., s_k, C)\subset M^{[n]}$ is also Lagrangian.
Under the natural map $\pi^{[n]}:\; M^{[n]}\arrow S^{[n]}=\C P^n$,
this subvariety is projected to a subset of $\Sym^nS=\C P^n$ consisting of
all $n$-tuples which contain $(s_1, ..., s_k)$.

\hfill

Another example is due to Richard Thomas 
(private communication). In the early versions of this
paper, we did not specify the behaviour of the restriction  $\pi\restrict Z$ 
outside of its regular values, and this example shows that it can be pretty wild.

Let $S$ be a compact complex 
torus, $\dim_\C S=n$, and $M= T^* S$ the total space of
its  cotangent bundle. Since $T^*S$ admits a natural trivialization, the 
manifold $M$ is equipped with a Lagrangian fibration $M \arrow \C^n$,
with the fibers obtained as translates of $S$.

Let $X\subset S$ be a complex submanifold, and $Z:=NS^\bot\subset M$
the total space of its conormal bundle. It is always Lagrangian,
and in many situations the projection of $Z$ to the base $\C^n$
is $n$-dimensional. To illustrate it, let us identify the base
$B=\C^n$ of $\pi$ with $T^*_s S$, for some $s\in S$.
A vector $v\in B$ belongs to $\pi(Z)$
if and only if $v\in T_xX^\bot$ for some $x\in X$,
where $T_xX^\bot= \{\zeta\in T^*_x S\ \ |\ \ \langle \zeta, T_xX\rangle =0\}.$
Then $\pi(Z)$ is a union of subspaces $T_xX^\bot$
parametrized by the family of $x\in X$. If, for example,
$X$ is a curve, and $T_x X$ is not constant, this
is a union of a non-constant family of hyperplanes,
hence it is Zariski dense in $B$.

Unless the tangent space $T_xX$ stays constant
as we vary $x\in X$, the image $\pi(Z)$ is $n$-dimensional,
and the corresponding fiber is 0-dimensional.
However, the central fiber $\pi^{-1}(0)$ is $X$, 
not a torus and of different dimension from the general fiber.
The proof of \ref{_subva_main_Theorem_} fails for the central fiber, because
$\pi(Z)$ is not smooth at $0$, and the identification
$T_z Z_x=\pi(T_z Z)^\bot \cap T^\pi_zM$ does not hold.

\hfill

{\bf Acknowledgments.}
Great many thanks to Richard Thomas for his comments, questions
and examples. We are indebted to Nigel Hitchin for his interest and inspiration,
and to Laura Schaposnik for her comments and expertise.
We express our gratitude to Stony Brook University and to the SCGP,
where this paper was prepared, and for their hospitality.
We thank the referees for their suggestions. 
Data sharing not applicable to this article as no datasets were generated 
or analysed during the current study.

{\small 
\noindent {\sc Ljudmila Kamenova\\
Department of Mathematics, 3-115 \\
Stony Brook University \\
Stony Brook, NY 11794-3651, USA,} \\
\tt kamenova@math.stonybrook.edu
\\

\noindent {\sc Misha Verbitsky\\
{\sc Instituto Nacional de Matem\'atica Pura e
              Aplicada (IMPA) \\ Estrada Dona Castorina, 110,
Jardim Bot\^anico, CEP 22460-320,\\ Rio de Janeiro, RJ - Brasil }\\
also:\\
{\sc Laboratory of Algebraic Geometry, HSE University,\\
Department of Mathematics, 9 Usacheva str., Moscow, Russia,}\\
\tt  verbit@mccme.ru}.
 }

\end{document}